\documentclass{article}
\usepackage[preprint,nonatbib]{neurips_2025}
\usepackage[T1]{fontenc}
\usepackage[utf8]{inputenc}
\usepackage{amsmath,amssymb,amsthm,mathtools}
\usepackage{booktabs}
\usepackage{microtype}
\usepackage{tikz}
\usepackage[hidelinks]{hyperref}

\newtheorem{theorem}{Theorem}[section]
\newtheorem{proposition}[theorem]{Proposition}
\newtheorem{corollary}[theorem]{Corollary}
\newtheorem{lemma}[theorem]{Lemma}
\theoremstyle{definition}
\newtheorem{definition}[theorem]{Definition}
\theoremstyle{remark}

\newcommand{\op}{\mathrm{op}}
\newcommand{\R}{\mathbb R}
\newcommand{\Sn}{\mathfrak S_n}
\newcommand{\norm}[1]{\left\lVert#1\right\rVert_{\op}}
\newcommand{\fall}[2]{(#1)_{#2}}
\newcommand{\rise}[2]{(#1)^{\overline{#2}}}
\newcommand{\ind}{\mathbf 1}
\title{A Graph-Theoretic Approach to Matrix AM--GM:\\Counterexamples and Sharp Exponential Bounds}
\author{Zehua Liu\\
  The University of Hong Kong\\
  \texttt{liuzehua@connect.hku.hk}}
\date{}

\begin{document}
\maketitle

\begin{abstract}
We develop a graph-theoretic approach to the expected-norm matrix
arithmetic--geometric mean inequality. By encoding matrix interactions with a
path-supported Gram matrix, we obtain a projection family whose product is
nonzero exactly when the transition graph of its index sequence is a path.
This gives an intuitive explanation for failure: repeated indices can create
cycles or branches that annihilate products, an effect absent for positive
scalars. We prove the universal lower bound $(n)_m/n^m$ for the ratio of the
with-replacement to the without-replacement mean and show that it is
asymptotically sharp on the exponential scale when $m=n$. Specifically, our
construction gives a ratio $\exp(-n+O(\sqrt n))$, matching the lower bound's
leading exponential rate over unrestricted dimensions. We also determine an
explicit overlap threshold, derive exact six-factor counterexamples, and show
that exponential failure can coexist with vanishing pairwise commutators.
\end{abstract}

\section{Introduction}
\label{sec:intro}

Given symmetric positive semidefinite matrices $A_1,\ldots,A_n$ and a
product length $1\le m\le n$, compare two averages of
$\|A_{i_1}\cdots A_{i_m}\|_{\op}$, where $\|\cdot\|_{\op}$ is the
Euclidean operator norm. Let $W_{n,m}$ be the average when the indices are
drawn independently and uniformly with replacement, and $D_{n,m}$ the average
when they are drawn uniformly without replacement. The expected-norm matrix
arithmetic--geometric mean (AM--GM) inequality asks whether $W_{n,m}\ge D_{n,m}$.
In both averages, the norm is taken before averaging; the full formulas appear
in Section~\ref{sec:problem}.

For nonnegative scalars and $m=n$, this comparison is the classical AM--GM
inequality: the two means are $(n^{-1}\sum_i a_i)^n$ and $\prod_i a_i$.
The matrix comparison can fail~\cite{desa2020}. We ask why repeated indices
can reverse the scalar comparison, and how small the ratio
$R_{n,m}=W_{n,m}/D_{n,m}$ can become when $D_{n,m}>0$.

Our construction isolates a geometric mechanism. Connect two distinct labels
whenever they occur consecutively in an index sequence; these edges form its
\emph{transition graph}. A sequence with distinct labels gives a path.
Repeating labels can preserve a path, close a cycle, or create a branch. We
construct projections whose product is nonzero exactly when this graph is a
path. For example, $(a,b,c,a)$ creates a triangle and gives a zero product,
whereas $(a,b,a,c)$ gives the path $b-a-c$ and survives. Some repetitions also
increase product norms: a consecutive repetition of a projection causes no
additional contraction. Thus the comparison depends on the balance between
repeated sequences that vanish and those that survive with increased weight.
The graph criterion makes both contributions countable.

The main quantitative result is a universal lower bound and a construction
matching its leading exponential rate. Every matrix family satisfies
\[
 R_{n,m}\ge\frac{n!}{(n-m)!\,n^m},
\]
by retaining the distinct-index terms. Our projections, in dimension $n\,n!$,
satisfy $R_{n,n}=\exp(-n+O(\sqrt n))$. Hence the lower bound has the optimal
leading exponential rate over unrestricted dimensions
(Theorem~\ref{thm:exponential} and Corollary~\ref{cor:optimal}). The same
counting framework gives exact six-factor counterexamples, an explicit
threshold as neighboring-vector inner products vary, and exponential failure
with vanishing pairwise commutators.

After related work and preliminaries, Section~\ref{sec:construction} proves
the graph criterion and Section~\ref{sec:bounds} establishes the main bounds.
Sections~\ref{sec:fixed}--\ref{sec:robustness} develop the concrete examples,
threshold, and robustness consequences.

\section{Related work}
\label{sec:related}

Recht and R\'e~\cite{recht2012} studied noncommutative AM--GM inequalities
motivated by sampling comparisons in randomized algorithms.
Israel, Krahmer, and Ward~\cite{israel2015} proved the expected-norm
inequality for products of up to three positive semidefinite matrices under
any unitarily invariant norm, including the operator norm.
De Sa~\cite{desa2020}, Section~2.3, gave an explicit expected-operator-norm
counterexample with six matrices and product length six, using
two-dimensional projections combined by direct sums over permutations.

Our construction also uses permutation direct sums. Its distinguishing feature
is a path-supported Gram matrix, which makes survival of a product equivalent
to a precise condition on its transition graph. This criterion explains the
role of repeated indices and yields quantitative bounds, including the optimal
leading exponential rate over unrestricted dimensions and an overlap threshold.
Our six-factor examples illustrate this framework; they do not improve the
matrix count in De Sa's example. Throughout, we study the average of product
norms, with the norm taken before averaging.

\section{Problem formulation and preliminaries}
\label{sec:problem}

\subsection{Matrix notation and the two sampling schemes}

A real matrix $A$ is \emph{symmetric positive semidefinite}, written
$A\succeq0$, if $A=A^\top$ and $x^\top Ax\ge0$ for every vector $x$.
The notation $A\preceq C$ means $C-A\succeq0$; $A\succ0$ means
$x^\top Ax>0$ for every nonzero $x$. We use the Euclidean operator norm,
\[
 \norm M=\sup_{\|x\|_2=1}\|Mx\|_2.
\]
The factors below are positive semidefinite, although their products need not
be symmetric. The operator norm is therefore defined by vector amplification,
rather than by an eigenvalue formula for the products.

Let $A_1,\ldots,A_n\in\R^{d\times d}$ be symmetric positive semidefinite
matrices, and let $1\le m\le n$. Write $\norm{\cdot}$ for the Euclidean
operator norm and define
\begin{align}
 W_{n,m}(A)&=\frac1{n^m}\sum_{i_1,\ldots,i_m=1}^n
       \norm{A_{i_1}\cdots A_{i_m}},\label{eq:W}\\
 D_{n,m}(A)&=\frac1{\fall nm}
       \sum_{\substack{i_1,\ldots,i_m\in\{1,\ldots,n\}\\\text{all distinct}}}
       \norm{A_{i_1}\cdots A_{i_m}},\label{eq:D}
\end{align}
where $\fall nm=n(n-1)\cdots(n-m+1)$. Whenever $D_{n,m}(A)>0$, set
\begin{equation}
 R_{n,m}(A)=\frac{W_{n,m}(A)}{D_{n,m}(A)}.
 \label{eq:ratio}
\end{equation}
The expected-norm inequality asserts that $W_{n,m}(A)\ge D_{n,m}(A)$, or
equivalently $R_{n,m}(A)\ge1$. Throughout the paper, the norm is taken
\emph{before} averaging. Taking the norm of an average instead gives a different
matrix AM--GM formulation.

Equivalently, let $I_1,\ldots,I_m$ be independent uniform labels in
$[n]=\{1,\ldots,n\}$, and let $\mathsf E$ be the event that they are all
distinct. Then
\[
 W_{n,m}(A)=\mathbb E\norm{A_{I_1}\cdots A_{I_m}},\qquad
 D_{n,m}(A)=\mathbb E\bigl[\norm{A_{I_1}\cdots A_{I_m}}\mid\mathsf E\bigr].
\]
Conditioning makes all ordered distinct tuples equally likely. This explains
the denominator $(n)_m$ in~\eqref{eq:D}; in particular, $(n)_n=n!$.
The two means are equal when $m=1$. We keep $n$ (the number of matrices),
$m$ (the product length), and $d$ (the matrix dimension) separate: some results
fix $m$ while $n$ grows, and others impose $m=n$.

As a scalar check, when $n=m=2$ and the matrices are nonnegative numbers
$a,b$, the two means are
\[
 W_{2,2}=(a^2+2ab+b^2)/4,\qquad D_{2,2}=ab,
 \qquad W_{2,2}-D_{2,2}=(a-b)^2/4\ge0.
\]
The numerator of the matrix comparison likewise includes repeated-index
products, but their norms are governed by relative vector geometry rather
than ordinary scalar multiplication.

\subsection{Projections, Gram matrices, and direct sums}
\label{sec:buildingblocks}

For a unit vector $u$, the matrix $P=uu^\top$ is the orthogonal projection
onto its span: $Px=\langle u,x\rangle u$. It satisfies $P^2=P=P^\top$,
$0\preceq P\preceq I$, and $\norm P=1$. Thus consecutive repetition of
the same projection adds no contraction. For two different unit vectors,
\[
 (uu^\top)(vv^\top)=\langle u,v\rangle uv^\top,
 \qquad \norm{(uu^\top)(vv^\top)}=|\langle u,v\rangle|.
\]
In particular, orthogonal vectors give a zero product. This identity is the
link between vector geometry and the graph construction.

The \emph{Gram matrix} of vectors $u_1,\ldots,u_n$ has entries
$G_{ij}=\langle u_i,u_j\rangle$. It is positive semidefinite because
$x^\top Gx=\|\sum_i x_i u_i\|_2^2\ge0$. Conversely, every real symmetric
positive semidefinite $G$ is a Gram matrix: its positive semidefinite square
root gives vectors $u_i=G^{1/2}e_i$, where $e_i$ is the $i$th standard basis
vector. Diagonal entries equal to one ensure that these vectors are unit
vectors. This lets us prescribe inner products before constructing vectors.

The \emph{direct sum} $M_1\oplus\cdots\oplus M_s$ is the block diagonal
matrix with diagonal blocks $M_1,\ldots,M_s$. Products with the same block
structure multiply block by block, and
\[
 \norm{\bigoplus_{j=1}^s M_j}=\max_{1\le j\le s}\norm{M_j}.
\]
Indeed, on a vector partitioned into blocks, the squared output norm is the
sum of the squared block output norms. This gives the upper bound by the
largest block norm, and a vector supported in a maximizing block gives
equality. In our construction, each block represents a different relabeling
of the same vectors. The direct sum permits a product to select its best
relabeling without changing the matrix family between products.

\section{From graphs to matrix products}
\label{sec:construction}

\subsection{The graph behind an index sequence}

The starting point is to record which distinct labels must interact in a
product. This separates the combinatorial pattern of repeated indices from
the matrix realization of that pattern.

We use finite undirected simple graphs: an edge joins two different vertices,
and there are no multiple edges. A \emph{path} on $k$ vertices has an ordering
$v_1,\ldots,v_k$ with exactly the edges $\{v_j,v_{j+1}\}$.
A \emph{cycle} is a closed chain of at least three different vertices;
a \emph{branch} here means a vertex with at least three neighbors.
An embedding into a path is an injective map of vertices sending every edge
to an edge of that path. No requirement is imposed on pairs that are not
edges. A connected simple graph embeds into a sufficiently long path precisely
when it has no cycles and no branches.

\begin{definition}
An index sequence $w=(i_1,\ldots,i_m)\in[n]^m$ is called a \emph{word}.
Let $r(w)=\#\{t\in\{1,\ldots,m-1\}:i_t\ne i_{t+1}\}$ be its
number of changes of label. Its \emph{transition graph} $\Gamma_w$ has
the labels occurring in $w$ as vertices and the distinct unordered pairs
$\{i_t,i_{t+1}\}$ with $i_t\ne i_{t+1}$ as edges. Repeated edges are ignored.
A single vertex is regarded as a path.
\end{definition}

We call a repeated label an \emph{index collision}. A \emph{hold} is a step
with $i_t=i_{t+1}$, so a word has $m-1-r(w)$ holds. A hold creates no graph
edge. Moving along an edge already traversed does not create a second edge,
but it still counts toward $r(w)$. The graph will decide whether a product
vanishes, whereas $r(w)$ will decide its nonzero weight. For example,
$(a,b,a,c)$ has three changes but just two edges, $\{a,b\}$ and $\{a,c\}$,
forming the path $b-a-c$.

Without replacement, this graph is always a path. With replacement, identifying
repeated labels may preserve a path, create a cycle, or produce a branch.
We choose a path as the graph of allowed interactions and realize each edge as
a nonzero inner product. All nonedges correspond to orthogonal vectors.
Taking a direct sum over labelings ensures that a word survives whenever its
transition graph admits \emph{any} embedding into the allowed path.

\subsection{Realizing the graph by projections}

Let $H_n$ denote the adjacency matrix of the path $1-2-\cdots-n$. For
$0<\rho\le1/2$, define
\begin{equation}
 G_n(\rho)=I_n+\rho H_n.
 \label{eq:gram}
\end{equation}
Here an adjacency matrix has entry one for neighboring vertices and zero
otherwise, including on the diagonal. Thus $G_n(\rho)$ has diagonal entries
one, neighboring entries $\rho$, and all other entries zero. We call $\rho$
the \emph{overlap parameter}: it is the inner product of neighboring vectors,
not an edge probability. The range $0<\rho\le1/2$ is a sufficient range for
positive semidefiniteness for every $n$. Indeed, for any $x\in\R^n$,
\[
 x^\top G_n(\rho)x
 \ge \sum_{i=1}^n x_i^2-\rho\sum_{i=1}^{n-1}(x_i^2+x_{i+1}^2)\ge0.
\]
Consequently, there are unit vectors $u_1,\ldots,u_n\in\R^n$ such that
\begin{equation}
 \langle u_i,u_j\rangle=
 \begin{cases}
  1,&i=j,\\
  \rho,&|i-j|=1,\\
  0,&|i-j|\ge2.
 \end{cases}
 \label{eq:inner}
\end{equation}
For example, take $u_i=G_n(\rho)^{1/2}e_i$. Let $P_i=u_iu_i^\top$.
In a product of these rank-one projections, consecutive labels contribute
their inner product, so a forbidden interaction makes the product zero.
To allow every relabeling of the graph simultaneously, form the direct sums
\begin{equation}
 B_i=\bigoplus_{\sigma\in\Sn}P_{\sigma(i)},\qquad i=1,\ldots,n,
 \label{eq:blocks}
\end{equation}
where $\Sn$ is the symmetric group on $\{1,\ldots,n\}$. Each $B_i$ is an
orthogonal projection of rank $n!$ acting on $\R^{n\,n!}$, and
$0\preceq B_i\preceq I$. At $\rho=1/2$, the alternative realization
$u_i=(e_i+e_{i+1})/\sqrt2$ in $\R^{n+1}$ gives an explicit coordinate
construction. We write $B=B^{(n,\rho)}$ and suppress parameters when clear.

For a concrete view of the block index, when $n=3$ the blocks of $B_1$ are
the six projections $P_{\sigma(1)}$, one for each permutation $\sigma$;
the corresponding blocks of $B_2$ and $B_3$ use the \emph{same} permutation.
Consequently a product has blocks
$P_{\sigma(i_1)}\cdots P_{\sigma(i_m)}$. A single block may annihilate a
word because of its particular labeling. The whole product vanishes only
when every permutation block annihilates it. This is why all permutations
are included.

\subsection{An exact graph criterion}

\begin{theorem}[Graph criterion for products]
\label{thm:graph}
For the projections in~\eqref{eq:blocks} and every word $w=(i_1,\ldots,i_m)$,
\begin{equation}
 \norm{B_{i_1}\cdots B_{i_m}}=
 \begin{cases}
   \rho^{r(w)},&\Gamma_w\text{ is a path},\\
   0,&\text{otherwise}.
 \end{cases}
 \label{eq:graph}
\end{equation}
In particular, for $1\le m\le n$,
\begin{equation}
 D_{n,m}(B)=\rho^{m-1}.
 \label{eq:denominator}
\end{equation}
\end{theorem}

\begin{proof}
Products of rank-one projections satisfy
\[
 P_{j_1}\cdots P_{j_m}
 =u_{j_1}\left(\prod_{t=1}^{m-1}\langle u_{j_t},u_{j_{t+1}}\rangle\right)u_{j_m}^\top.
\]
The endpoint vectors are unit vectors, and the operator norm of a direct sum is
the maximum of its block norms. Hence
\begin{equation}
 \norm{B_{i_1}\cdots B_{i_m}}
 =\max_{\sigma\in\Sn}\prod_{t=1}^{m-1}
   |\langle u_{\sigma(i_t)},u_{\sigma(i_{t+1})}\rangle|.
 \label{eq:max}
\end{equation}
By~\eqref{eq:inner}, a block is nonzero precisely when $\sigma$ sends every
edge of $\Gamma_w$ to an edge of the ambient path. The transition graph is
connected on its vertex set: segments of the word provide walks between
its labels. A connected simple graph admits an injective
edge-preserving map into a path if and only if it is itself a path. Any such map
extends to a permutation of all $n$ labels. When it exists, each change of label
contributes $\rho$ and each hold contributes $1$, proving~\eqref{eq:graph}.
For a word with distinct labels, the transition graph is always a path and
$r(w)=m-1$, proving~\eqref{eq:denominator}.
\end{proof}

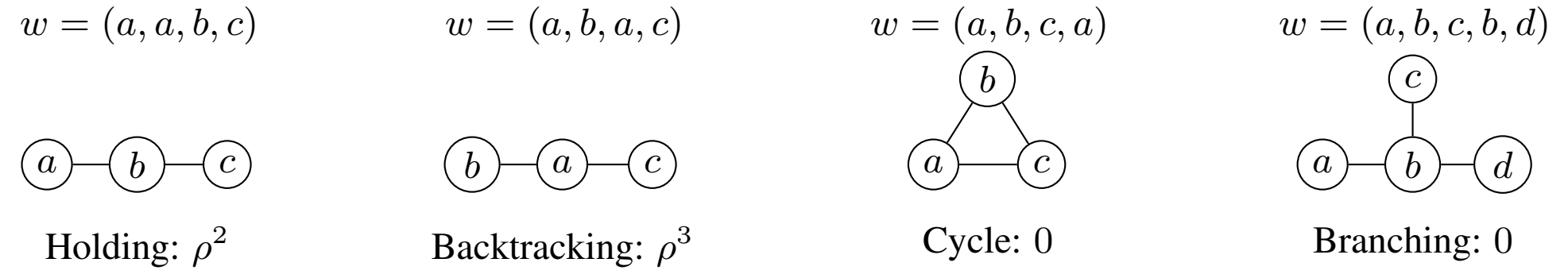
\begin{figure}[t]
\centering
\begin{minipage}[t]{.24\textwidth}
\centering
\small $w=(a,a,b,c)$\\[4pt]
\begin{tikzpicture}[baseline=0,scale=.75,every node/.style={circle,draw,inner sep=2pt,font=\small}]
\path[use as bounding box] (-.3,-.55) rectangle (2.3,1.1);
\node (a) at (0,0) {$a$}; \node (b) at (1,0) {$b$}; \node (c) at (2,0) {$c$};
\draw (a)--(b)--(c);
\end{tikzpicture}\\[2pt]
Holding: $\rho^2$
\end{minipage}\hfill
\begin{minipage}[t]{.24\textwidth}
\centering
\small $w=(a,b,a,c)$\\[4pt]
\begin{tikzpicture}[baseline=0,scale=.75,every node/.style={circle,draw,inner sep=2pt,font=\small}]
\path[use as bounding box] (-.3,-.55) rectangle (2.3,1.1);
\node (b) at (0,0) {$b$}; \node (a) at (1,0) {$a$}; \node (c) at (2,0) {$c$};
\draw (b)--(a)--(c);
\end{tikzpicture}\\[2pt]
Backtracking: $\rho^3$
\end{minipage}\hfill
\begin{minipage}[t]{.24\textwidth}
\centering
\small $w=(a,b,c,a)$\\[4pt]
\begin{tikzpicture}[baseline=0,scale=.75,every node/.style={circle,draw,inner sep=2pt,font=\small}]
\path[use as bounding box] (-.3,-.55) rectangle (2.3,1.1);
\node (a) at (.4,0) {$a$}; \node (b) at (1,.95) {$b$}; \node (c) at (1.6,0) {$c$};
\draw (a)--(b)--(c)--(a);
\end{tikzpicture}\\[2pt]
Cycle: $0$
\end{minipage}\hfill
\begin{minipage}[t]{.24\textwidth}
\centering
\small $w=(a,b,c,b,d)$\\[4pt]
\begin{tikzpicture}[baseline=0,scale=.75,every node/.style={circle,draw,inner sep=2pt,font=\small}]
\path[use as bounding box] (-.3,-.55) rectangle (2.3,1.1);
\node (a) at (0,0) {$a$}; \node (b) at (1,0) {$b$}; \node (c) at (1,.95) {$c$}; \node (d) at (2,0) {$d$};
\draw (a)--(b)--(c); \draw (b)--(d);
\end{tikzpicture}\\[2pt]
Branching: $0$
\end{minipage}
\caption{Index collisions and their algebraic consequences. The diagrams show
transition graphs, and the values below are the corresponding product norms.
Holding and backtracking preserve paths. Cycles and branches force zero
products under every labeling.}
\label{fig:graphs}
\end{figure}

Figure~\ref{fig:graphs} illustrates the criterion. A consecutive repetition
contributes a hold of weight $1$, whereas moving along an allowed edge
contributes $\rho$. Backtracking can therefore preserve a product. A cycle or
a degree-three vertex, however, cannot be embedded into the interaction path,
so the corresponding product vanishes in every block.

For example, take $n\ge3$. The word $(1,2,3,1)$ requires three pairwise
interactions and therefore a triangle; no block can realize all three edges,
so its product norm is zero. By contrast, $(1,2,1,3)$ requires only the path
$2-1-3$. A block placing label $1$ between labels $2$ and $3$ survives and
has norm $\rho^3$. Finally, $(1,1,2,3)$ has a hold and only two changes,
giving norm $\rho^2$. Thus both the graph and the number of changes matter.

This explains how the scalar intuition can fail. Although some repeated words
have greater weight than words with distinct labels, many repeated words have
zero weight. To compare the two means, one must count both effects. The
universal lower bound comes from the distinct-label words, and the graph
criterion lets us show that the surviving repeated words add only a
subexponential factor when $m=n$ and $\rho$ is fixed.

\section{A universal lower bound and its exponential sharpness}
\label{sec:bounds}

Once $R_{n,m}\ge1$ fails, a natural replacement is a smaller universal lower
bound. The sampling interpretation gives such a bound for every matrix family.

\begin{theorem}[Universal lower bound]
\label{thm:universal}
Let $1\le m\le n$. For any matrices $A_1,\ldots,A_n$ with
$D_{n,m}(A)>0$,
\begin{equation}
 R_{n,m}(A)\ge p_{n,m}:=\frac{\fall nm}{n^m}.
 \label{eq:universal}
\end{equation}
In particular, this bound holds for every positive semidefinite matrix family,
independently of dimension.
\end{theorem}

\begin{proof}
The with-replacement sum contains every without-replacement term, and all
remaining norms are nonnegative. Consequently,
$n^mW_{n,m}(A)\ge\fall nm D_{n,m}(A)$.
\end{proof}

The coefficient $p_{n,m}$ is exactly the probability that $m$ independent
uniform draws from $n$ labels are all distinct. We next bound the extra
contribution of repeated words for the graph construction and prove that,
when $m=n$, this lower bound has the optimal leading exponential rate.

\subsection{Counting the surviving words}

The conditioning interpretation makes the remaining task precise. For $m\ge2$,
let $Q_{n,m}$ be the conditional mean product norm given at least one collision.
Then
\[
 R_{n,m}=p_{n,m}+(1-p_{n,m})\frac{Q_{n,m}}{D_{n,m}}.
\]
Thus $R_{n,m}<1$ exactly when $Q_{n,m}<D_{n,m}$. To approach the universal
lower bound, we must control the total contribution of repeated words relative
to the distinct-word contribution. The graph criterion enables us to do this
without constructing the large block matrices numerically.

Write $\rise a\ell=a(a+1)\cdots(a+\ell-1)$ for the rising factorial, with
$\rise a0=1$.
In the count below, $k$ denotes the number of distinct labels in a word and
$\ell=m-k$ the number of occurrences beyond the first occurrences. This
$\ell$ includes revisits after moving away; it is not just the number of holds.

\begin{theorem}[A uniform ratio bound]
\label{thm:count}
Let $C_\rho=2+\rho^{-1}$. For $1\le m\le n$, the path-Gram family satisfies
\begin{equation}
 p_{n,m}\le R_{n,m}(B)\le
 p_{n,m}\sum_{\ell=0}^{m-1}\binom{m-1}{\ell}
       \frac{C_\rho^\ell}{\rise{n-m+1}\ell}.
 \label{eq:count}
\end{equation}
\end{theorem}

\begin{proof}
The lower bound is~\eqref{eq:universal}. For the upper bound, classify words
by their number $k$ of distinct labels. The first-occurrence positions after
the initial position can be chosen in $\binom{m-1}{k-1}$ ways, and assigning
labels in order of first appearance gives $\fall nk$ choices. Each new label
contributes a factor $\rho$.
For example, $(a,a,b,a,c)$ first introduces labels at positions $1,3,5$.
Once these positions and the ordered labels $(a,b,c)$ are fixed, the other
positions must use labels already introduced. This explains the position
count $\binom{m-1}{k-1}$ rather than $\binom mk$.

If the final transition graph is a path, then every prefix also has a path as
its transition graph. At a repeated-label step, one can hold at the current
vertex or move to at most two existing neighbors. The sum of their weights is
at most $1+2\rho$. Moving to a previously visited vertex not already adjacent
to the current vertex creates a cycle in the connected prefix graph, and is
therefore forbidden. Introducing a new label can also create a branch; ignoring
that restriction only increases the bound. It follows that the sum of the norms
of all nonzero words using exactly $k$ labels is at most
\begin{equation}
 \binom{m-1}{k-1}\fall nk\,\rho^{k-1}(1+2\rho)^{m-k}.
 \label{eq:kbound}
\end{equation}
There are $k-1$ new-label steps and $m-k$ repeated-label steps. Dividing
the weight bound by $D_{n,m}=\rho^{m-1}$ gives
\[
 \rho^{k-1}(1+2\rho)^{m-k}/\rho^{m-1}
 =\bigl((1+2\rho)/\rho\bigr)^{m-k}=C_\rho^{m-k}.
\]
Sum over $k$, divide by $n^m$, and set $\ell=m-k$. The identity
$\fall n{m-\ell}=\fall nm/\rise{n-m+1}\ell$ gives~\eqref{eq:count}.
\end{proof}

\subsection{Asymptotic attainment of the lower bound}

We now set $m=n$ and keep $\rho$ fixed while $n$ grows. A sequence of positive
factors $h_n$ is called \emph{subexponential} when $\log h_n=o(n)$.
Multiplying by such a factor does not change a limit of $n^{-1}\log R_n$.
We use $O_\rho(\sqrt n)$ for a quantity bounded in absolute value by a
constant depending on $\rho$ times $\sqrt n$, and $o(n)$ for a quantity
whose ratio to $n$ tends to zero. The next theorem allows a subexponential
gap from the lower bound, which is the precise meaning of exponential
sharpness in this paper.

\begin{theorem}[Matching the lower bound on the exponential scale]
\label{thm:exponential}
Fix $0<\rho\le1/2$. For products of length $n$,
\begin{equation}
 \frac{n!}{n^n}\le R_{n,n}(B)\le
 \frac{n!}{n^n}\exp\!\left(2\sqrt{C_\rho(n-1)}\right).
 \label{eq:expbound}
\end{equation}
Consequently, $R_{n,n}(B)=\exp(-n+O_\rho(\sqrt n))$ as $n\to\infty$, and
\begin{equation}
 \lim_{n\to\infty}\frac1n
     \log\frac{R_{n,n}(B)}{p_{n,n}}=0.
 \label{eq:sharpness}
\end{equation}
\end{theorem}

\begin{proof}
Taking $m=n$ in Theorem~\ref{thm:count} gives
\[
 R_{n,n}(B)\le\frac{n!}{n^n}
 \sum_{\ell=0}^{n-1}\binom{n-1}{\ell}\frac{C_\rho^\ell}{\ell!}.
\]
Use $\binom{n-1}{\ell}\le(n-1)^\ell/\ell!$ and, for $x\ge0$,
\[
 \sum_{\ell\ge0}\frac{x^\ell}{(\ell!)^2}
 \le e^{\sqrt x}e^{\sqrt x}=e^{2\sqrt x}.
\]
The inequality follows by retaining only the diagonal terms in the product of
the two exponential series. This proves the upper bound. Stirling's formula
$\log(n!/n^n)=-n+O(\log n)$, together with the lower bound, proves the claimed
asymptotic behavior.
\end{proof}

To state the worst-case consequence, define the optimal universal coefficient
with unrestricted dimension by
\begin{equation}
 \mathcal C_n=\inf_{\substack{d\ge1,\ A_1,\ldots,A_n\succeq0\\D_{n,n}(A)>0}}
                   R_{n,n}(A).
 \label{eq:bestconstant}
\end{equation}

\begin{corollary}[Sharpness over all dimensions]
\label{cor:optimal}
For every $n$,
\begin{equation}
 \frac{n!}{n^n}\le\mathcal C_n\le
 \frac{n!}{n^n}e^{4\sqrt{n-1}},
 \qquad
 \lim_{n\to\infty}\frac{\log\mathcal C_n}{n}=-1.
 \label{eq:optimal}
\end{equation}
\end{corollary}

\begin{proof}
The lower bound holds for all matrix families. The upper bound follows from
Theorem~\ref{thm:exponential} with $\rho=1/2$.
\end{proof}

Thus the universal lower bound and the worst-case ratio have the same leading
exponential rate. The result does not assert the stronger equivalence
$R_{n,n}(B)/p_{n,n}\to1$; the multiplicative gap in~\eqref{eq:expbound} is
subexponential in $n$.

In particular, no positive constant $c$ independent of $n$ can make
$W_{n,n}\ge cD_{n,n}$ hold for all positive semidefinite matrix families.
The separation can also be normalized while keeping individual matrix norms
uniformly bounded. At $\rho=1/2$, set
$\widetilde B_i=2^{(n-1)/n}B_i$. Homogeneity gives
\begin{equation}
 D_{n,n}(\widetilde B)=1,\qquad
 W_{n,n}(\widetilde B)=\exp(-n+O(\sqrt n)),\qquad
 \norm{\widetilde B_i}<2.
 \label{eq:normalization}
\end{equation}

\section{Counterexamples at fixed product length}
\label{sec:fixed}

The main bound takes a complete pass, $m=n\to\infty$. We now hold the
product length $m$ and overlap $\rho$ fixed and increase only the number $n$
of matrices. In this regime, almost all sampled words have distinct labels,
so the ratio tends to one. Its first correction reveals whether it approaches
one from above or below. Patterns with just one repeated pair determine this
correction; words with two or more repeated occurrences enter only at order
$n^{-2}$.

\begin{theorem}[Fixed-length expansion]
\label{thm:fixed}
Fix $m\ge4$ and $0<\rho\le1/2$. As $n\to\infty$,
\begin{equation}
 R_{n,m}(B)=1+
 \frac{(m-1)/\rho+2-\binom m2}{n}+O_{m,\rho}(n^{-2}).
 \label{eq:fixed}
\end{equation}
Thus, if $(m-1)/\rho+2<\binom m2$, then $R_{n,m}(B)<1$ for all
sufficiently large $n$. In particular, at $\rho=1/2$, every fixed $m\ge6$
gives counterexamples for arbitrarily large numbers of matrices.
\end{theorem}

\begin{proof}
Words with distinct labels contribute
$p_{n,m}=1-\binom m2/n+O_m(n^{-2})$. A word using exactly $m-1$ distinct
labels has one repeated pair of positions. If these positions are adjacent,
the transition graph is a path. There are $m-1$ choices of positions and
$\fall n{m-1}$ label assignments for each; the normalized weight is
$\rho^{-1}$.

Suppose instead that the repeated positions are $a<b$ with $b-a\ge2$.
If $b-a\ge3$, the intervening segment creates a cycle. If $b-a=2$ and both
$a>1$ and $b<m$, the repeated vertex has three distinct neighbors, creating
a branch. Only the pairs $(a,b)=(1,3)$ and $(m-2,m)$ preserve a path.
Each has $\fall n{m-1}$ label assignments and normalized weight $1$.

Words using at most $m-2$ labels number $O_m(n^{m-2})$. Their normalized
weights are uniformly bounded for fixed $m,\rho$, so their contribution is
$O_{m,\rho}(n^{-2})$. Combining these terms gives
\[
 R_{n,m}=p_{n,m}
  +\frac{\fall n{m-1}}{n^m}\bigl((m-1)/\rho+2\bigr)
  +O_{m,\rho}(n^{-2}),
\]
which proves~\eqref{eq:fixed}. At $\rho=1/2$, the first-order coefficient
is $m(5-m)/2$.
\end{proof}

\subsection{An exact weighted-walk formula}

The graph criterion also permits exact computation without forming the large
block-diagonal matrices. Set $q=\rho^{-1}$ and, for a given $m$, define
\begin{equation}
 F_k=\ind^\top(qI_k+H_k)^{m-1}\ind,\qquad F_0=0,
 \label{eq:walks}
\end{equation}
where $\ind$ denotes the all-ones vector of the appropriate dimension.

The matrix $qI_k+H_k$ is a weighted transition matrix on the $k$-vertex
path: staying at a vertex has weight $q$, and moving to a neighbor has weight
one. Each entry of its $(m-1)$st power sums the weights of all walks between
the corresponding endpoints, with the weight of a walk defined as the product
of its step weights. Multiplication by the all-ones vectors sums over starting
and ending vertices. Thus $F_k$ counts weighted walks, rather than their
probabilities. For example, at $m=3$,
\[
 F_1=q^2,\qquad F_2=2(q+1)^2,\qquad F_2-2F_1=4q+2.
\]
The last expression counts the two-vertex walks that use both vertices:
four walks with one hold and two walks with no holds.

\begin{proposition}[Exact enumeration]
\label{prop:exact}
For $1\le m\le n$,
\begin{equation}
 n^mR_{n,m}(B)=nq^{m-1}
 +\frac12\sum_{k=2}^m\fall nk\,(F_k-2F_{k-1}+F_{k-2}).
 \label{eq:exact}
\end{equation}
\end{proposition}

\begin{proof}
After division by $D_{n,m}=\rho^{m-1}$, each hold has weight $q$ and each
move has weight $1$. The quantity $F_k$ counts all weighted walks of $m-1$
steps on a path with $k$ vertices. For $k\ge2$, inclusion--exclusion shows
that walks visiting both endpoints have total weight
$F_k-2F_{k-1}+F_{k-2}$: subtract the walks missing the left endpoint and
those missing the right endpoint, then add back those missing both. The
remaining graphs have $k-1$, $k-1$, and $k-2$ vertices, respectively.
Such walks have the entire path as their transition
graph. There are $\fall nk/2$ undirected labeled paths with $k$ vertices:
an ordered list of $k$ distinct labels specifies a path, and exactly its two
opposite orders specify the same undirected path when $k\ge2$.
Single-vertex words contribute $nq^{m-1}$, completing the count.
\end{proof}

\begin{corollary}[Explicit six-factor counterexamples]
\label{cor:six}
At $m=6$ and $\rho=1/2$,
\begin{align}
 n^6R_{n,6}(B)&=32\fall n1+211\fall n2+206\fall n3
                   +66\fall n4+12\fall n5+\fall n6,\label{eq:sixfall}\\
 R_{n,6}(B)&=1-\frac3n+\frac{31}{n^2}+\frac5{n^3}-\frac7{n^4}+\frac5{n^5}.
 \label{eq:sixpoly}
\end{align}
For every integer $n\ge11$, this ratio is strictly less than $1$.
\end{corollary}

\begin{proof}
Setting $q=2$ in~\eqref{eq:walks} gives
\[
 (F_1,F_2,F_3,F_4,F_5,F_6)=(32,486,1352,2350,3372,4396).
\]
Substitution into~\eqref{eq:exact} proves~\eqref{eq:sixfall}; expanding the
falling factorials gives~\eqref{eq:sixpoly}. For $n\ge11$,
\[
 n^5(R_{n,6}-1)=n^3(31-3n)+5n^2-7n+5
 \le-2n^3+5n^2-7n+5<0.
\]
\end{proof}

For example, $D_{11,6}=1/32$ and
$W_{11,6}/D_{11,6}=158922/161051<1$. Table~\ref{tab:six} lists values
obtained from the exact formula. The path-Gram family has ratio greater than
$1$ when $n=m=6$.

\begin{table}[ht]
\centering
\caption{Exact-formula evaluations for six-factor products at $\rho=1/2$.}
\label{tab:six}
\begin{tabular}{rrrrrr}
\toprule
$n$&6&10&11&12&20\\
\midrule
$R_{n,6}$&1.379501&1.014350&0.986781&0.967854&0.928083\\
\bottomrule
\end{tabular}
\end{table}

\section{A critical overlap scale}
\label{sec:transition}

When $\rho$ decreases with $n$, holding at the same label avoids increasingly
costly transitions. This effect competes with the annihilation of products
whose transition graphs are not paths. We now identify the threshold at the
scale $\rho\asymp n^{-1}$. Write
$R_n(\rho)=R_{n,n}(B^{(n,\rho)})$ and define
\begin{equation}
 T_n(\rho)=\frac{n!}{n^n}\sum_{\ell=0}^{n-1}
                    \binom{n-1}{\ell}\frac{\rho^{-\ell}}{\ell!}.
 \label{eq:T}
\end{equation}

Unlike the fixed-overlap regime of Section~\ref{sec:bounds}, here $m=n$ and
$\rho$ may change with $n$. The notation $\rho\asymp n^{-1}$ means that
$n\rho$ stays between positive constants. To see why this scale matters,
normalize each product norm by $D_{n,n}=\rho^{n-1}$. A surviving change of
label then has weight one, while each hold has relative weight $\rho^{-1}$.
Decreasing $\rho$ therefore favors words with many holds. The quantity
$T_n(\rho)$ retains a simple subset of such words, in which each label is
used in one consecutive run. The following lemma shows that this subset
already determines the exponential rate at $\rho=\kappa/n$.

\begin{lemma}
\label{lem:sandwich}
For $0<\rho\le1/2$,
\begin{equation}
 T_n(\rho)\le R_n(\rho)\le(1+2\rho)^{n-1}T_n(\rho).
 \label{eq:sandwich}
\end{equation}
\end{lemma}

\begin{proof}
Restrict attention to words in which every label occupies one contiguous run
and is never revisited. With $k$ runs, there are $\binom{n-1}{k-1}$ choices of
run lengths and $\fall nk$ assignments of distinct labels. Each product has
norm $\rho^{k-1}$. Their contribution to the ratio is
\[
 \frac1{n^n}\sum_{k=1}^n\binom{n-1}{k-1}\fall nk\,\rho^{k-n}=T_n(\rho),
\]
proving the lower bound. The upper bound~\eqref{eq:kbound} differs termwise
by a factor $(1+2\rho)^{n-k}$, which is at most $(1+2\rho)^{n-1}$.
\end{proof}

\begin{theorem}[Exact exponential rate and critical constant]
\label{thm:transition}
Fix $\kappa>0$ and set $\rho=\kappa/n$, with $n$ sufficiently large that
$\rho\le1/2$. Define
\begin{equation}
 t_\kappa=\frac{2}{1+\sqrt{1+4\kappa}},\qquad
 g(\kappa)=t_\kappa-1-\log(1-t_\kappa).
 \label{eq:g}
\end{equation}
Then
\begin{equation}
 \lim_{n\to\infty}\frac1n\log R_n(\kappa/n)=g(\kappa).
 \label{eq:rate}
\end{equation}
Let $\Omega$ be the unique positive solution of $\Omega e^\Omega=1$, and set
\begin{equation}
 \kappa_* =\frac{\Omega}{(1-\Omega)^2}
            =3.026944729\ldots.
 \label{eq:critical}
\end{equation}
If $0<\kappa<\kappa_*$, the ratio $R_n(\kappa/n)$ grows exponentially and
the expected-norm inequality eventually holds for this family. If
$\kappa>\kappa_*$, the ratio decays exponentially and the inequality eventually
fails strictly.
\end{theorem}

\begin{proof}
Since $(1+2\kappa/n)^{n-1}\le e^{2\kappa}$, Lemma~\ref{lem:sandwich}
shows that $R_n(\kappa/n)$ and $T_n(\kappa/n)$ have the same exponential
rate. For the $\ell$th summand in~\eqref{eq:T}, set $t=\ell/n$. Stirling's
formula gives the uniform expansion
\begin{align}
 \frac1n\log\left[\frac{n!}{n^n}\binom{n-1}{\ell}
                      \frac{(n/\kappa)^\ell}{\ell!}\right]
 &=\phi_\kappa(t)+O_\kappa\!\left(\frac{\log n}{n}\right),\label{eq:stirling}\\
 \phi_\kappa(t)&=-1-2t\log t-(1-t)\log(1-t)+t(1-\log\kappa),\nonumber
\end{align}
with $0\log0=0$. Uniformity follows from
$\log(j!)=j\log j-j+O(\log(n+1))$ for $0\le j\le n$, together with
$\binom{n-1}{\ell}=(1-\ell/n)\binom n\ell$. The logarithm of the last
factor lies between $-\log n$ and $0$ for $0\le\ell\le n-1$.
There are $n$ summands, so their sum lies between the largest summand and $n$
times that summand. Therefore
\begin{equation}
 \lim_{n\to\infty}\frac1n\log T_n(\kappa/n)
      =\max_{0\le t\le1}\phi_\kappa(t).
 \label{eq:laplace}
\end{equation}
Here $t$ is the proportion of repeated occurrences in the run-restricted
count, not the proportion of distinct labels. Equation~\eqref{eq:laplace}
reduces the asymptotic sum to a one-variable maximization; the factor $n$
between the sum and its largest summand disappears after taking $n^{-1}\log$.

On $(0,1)$,
\[
 \phi_\kappa'(t)=\log(1-t)-2\log t-\log\kappa,
 \qquad \phi_\kappa''(t)=-\frac1{1-t}-\frac2t<0.
\]
The unique maximizer solves $1-t=\kappa t^2$, hence equals $t_\kappa$.
Substituting gives $\phi_\kappa(t_\kappa)=g(\kappa)$ and proves~\eqref{eq:rate}.

Now put $s=1-t_\kappa$. Then $\kappa=s/(1-s)^2$ and
$g(\kappa)=-s-\log s$. On $(0,1)$, the former expression is strictly
increasing in $s$ and the latter strictly decreasing. The unique zero satisfies
$s e^s=1$, giving~\eqref{eq:critical} and the signs on either side.
\end{proof}

For instance, $g(2)=\log2-1/2=0.1931471806\ldots$, whereas
$g(4)=-0.1146788737\ldots$. At $\kappa=\kappa_*$, the theorem establishes
only a zero exponential rate; it does not determine the direction of the
inequality from the subexponential terms.

\begin{corollary}[Finite-size thresholds]
\label{cor:finite}
For all sufficiently large $n$, there is a unique $\rho_n\in(0,1/2)$ such
that $R_n(\rho_n)=1$, and
\begin{equation}
 n\rho_n\longrightarrow\kappa_*.
 \label{eq:finite}
\end{equation}
\end{corollary}

\begin{proof}
By Theorem~\ref{thm:graph}, $R_n(\rho)$ is a polynomial in $\rho^{-1}$
with nonnegative coefficients. For $n\ge2$, constant words contribute
$n^{1-n}\rho^{-(n-1)}$, so $R_n$ is strictly decreasing and tends to infinity
as $\rho\downarrow0$. Theorem~\ref{thm:exponential} gives $R_n(1/2)<1$
for sufficiently large $n$, proving existence and uniqueness of the root.
For any $0<\kappa_-<\kappa_*<\kappa_+$,
Theorem~\ref{thm:transition} implies
$R_n(\kappa_-/n)>1>R_n(\kappa_+/n)$ for sufficiently large $n$.
Monotonicity yields $\kappa_-<n\rho_n<\kappa_+$, proving the limit.
\end{proof}

\section{Approximate commutativity and positive-definite perturbations}
\label{sec:robustness}

Two further questions concern how much algebraic separation is needed for
failure. Must the matrices strongly fail to commute? Must they have zero
eigenvalues? The next two statements answer these separately. For matrices
$X,Y$, their \emph{commutator} is $XY-YX$; it vanishes exactly when they
commute. Positive definiteness, in contrast, excludes zero eigenvalues.

\begin{proposition}[Small commutators and exponential failure]
\label{prop:commutator}
For $n\ge2$, the path-Gram family satisfies
\begin{equation}
 \max_{i\ne j}\norm{B_iB_j-B_jB_i}=\rho\sqrt{1-\rho^2}.
 \label{eq:commutator}
\end{equation}
At $\rho=4/n$, these commutator norms tend to zero, whereas
\begin{equation}
 R_n(4/n)=\exp\!\left((-0.1146788737\ldots)n+o(n)\right).
 \label{eq:almost}
\end{equation}
\end{proposition}

\begin{proof}
For unit vectors $u,v$, let $c=|\langle u,v\rangle|$. A calculation on
their span gives
$\norm{uu^\top vv^\top-vv^\top uu^\top}=c\sqrt{1-c^2}$.
To see this for $c<1$, choose an orthonormal basis with
$u=e_1$ and, after changing the sign of $v$ if necessary,
$v=ce_1+\sqrt{1-c^2}e_2$. The commutator on this span is
\[
 \begin{pmatrix}0&c\sqrt{1-c^2}\\-c\sqrt{1-c^2}&0\end{pmatrix},
\]
whose operator norm is $c\sqrt{1-c^2}$; it is zero off the span.

In each direct-sum block, distinct vectors have inner product $0$ or $\rho$.
For every pair of distinct labels, some permutation sends them to neighboring
vertices. Taking the maximum block norm proves~\eqref{eq:commutator}.
Equation~\eqref{eq:almost} follows from Theorem~\ref{thm:transition}.
\end{proof}

Here the dimension, number of matrices, and product length all grow. Thus,
vanishing pairwise commutator norms alone do not ensure the expected-norm
inequality in this growing-length regime.

Finally, strict failure persists after a sufficiently small perturbation toward
the identity. The tolerance below depends on the finite-instance gap; it need
not stay bounded away from zero as $n$ grows. We therefore obtain
positive-definite counterexamples without asserting a uniform lower bound
on their smallest eigenvalues across the asymptotic family.

\begin{proposition}[Stability under positive-definite perturbations]
\label{prop:perturbation}
Suppose the path-Gram family at some finite $n,m,\rho$ satisfies
$\delta=D_{n,m}(B)-W_{n,m}(B)>0$. Define
$C_i(\varepsilon)=(1-\varepsilon)B_i+\varepsilon I$. If
$0<\varepsilon<\min\{1,\delta/(2m)\}$, then
$0\prec C_i(\varepsilon)\preceq I$ and $W_{n,m}(C)<D_{n,m}(C)$.
\end{proposition}

\begin{proof}
Both families consist of contractions and $\norm{C_i-B_i}\le\varepsilon$.
A telescoping expansion gives
\[
 \norm{C_{i_1}\cdots C_{i_m}-B_{i_1}\cdots B_{i_m}}\le m\varepsilon.
\]
By the reverse triangle inequality, each of $W_{n,m}$ and $D_{n,m}$ changes
by at most $m\varepsilon$. Hence
$D_{n,m}(C)-W_{n,m}(C)\ge\delta-2m\varepsilon>0$.
\end{proof}

\section{Discussion}

The graph-theoretic viewpoint makes the failure of matrix AM--GM explicit:
index collisions can change a path into a cycle or a branching graph and
thereby annihilate a product. This mechanism is invisible to a comparison of
individual matrix spectra. All matrices in the construction remain projections
of the same rank, with spectra, traces, and norms independent of $\rho$.
Their relative geometry controls the inequality.

The same viewpoint answers a quantitative question prompted by failure: the
collision-free probability $(n)_m/n^m$ is a universal lower bound, and for
$m=n$ it is attained at the leading exponential order by the graph
construction. Thus no universally larger exponential rate is possible when
dimension is unrestricted.

Open questions include the exact optimal coefficient and its subleading
asymptotics, its dependence on fixed dimension and product length, and the
subexponential behavior at $\kappa=\kappa_*$. Our results allow dimension to
grow and concern the operator norm. For a fixed Schatten $p$-norm, a direct sum
does not select the largest block; its enumeration changes, and the same
threshold does not follow from this argument.

\section*{Acknowledgments}

We acknowledge the use of OpenAI's GPT models to assist with developing and
checking the proofs presented in this paper.

\bibliographystyle{plainurl}
\bibliography{references}
\end{document}